\documentclass[2pt]{article}
\usepackage{graphicx}
\usepackage{amsmath}
\usepackage{mathtools}
\usepackage{mathrsfs}
\usepackage[hmarginratio=1:1,top=32mm,columnsep=20pt]{geometry}
\usepackage{tocloft}
\usepackage{dsfont}
\usepackage{tikz}
\usepackage{tabularx}
\usepackage{multicol}
\usetikzlibrary{arrows,shapes,positioning}
\usepackage{bbm}
\usepackage{bbold}
\usepackage{titlesec}
\usepackage{hyperref}
\hypersetup{
    colorlinks=true,
    linkcolor=blue,
    filecolor=magenta,
    urlcolor=cyan,
}
\usepackage{amsmath,amsfonts,amssymb,amsthm,epsfig,epstopdf,titling,url,array}

\theoremstyle{plain}
\newtheorem{thm}{Theorem}[section]

\newtheorem{prop}[thm]{Proposition}

\newtheorem{defn}[thm]{Definition}

\newtheorem{exmp}[thm]{Example}
\newtheorem{rem}[thm]{Remark}
\newtheorem{note}[thm]{Note}

\usepackage{amssymb}
\usepackage{authblk}
\usepackage{caption}
\usepackage{subcaption}

\title{\bf  The weak rate of  convergence for the Euler-Maruyama approximation  of one-dimensional stochastic differential equations involving the local times of the unknown process}
\author{ M.Benabdallah\thanks{M.Benabdallah : e-mail: bmohsine@gmail.com} }
\author{ K.Hiderah\thanks{K.Hiderah: e-mail: kamal200968@yahoo.com}  }
\affil{Department of Mathematics, faculty of science, University of Ibn Tofail, Kenitra, Morocco}
\date{}
\begin{document}
\maketitle
\begin{abstract}
In this paper, we consider the weak convergence  of the Euler-Maruyama approximation for  one dimensional stochastic differential equations  involving the local times of the unknown process.
We use a transformation in order to remove the local time $L_{t}^{a}$ from the stochastic differential equations of type:
\begin{equation*}
X_{t}=X_{0}+\int_{0}^{t}\varphi(X_{s})dB_{s}+\int_{\mathds{R}}\nu(da)L_{t}^{a}
\end{equation*}
 where $B$ is a one-dimensional Brownian motion, $\varphi:\mathds{R}\rightarrow \mathds{R}$ is a bounded measurable function, and $\nu$ is a
 bounded measure on $\mathds{R}$ and we provide the approximation of Euler-maruyama for the stochastic differential equations without local time.
 After that, we conclude the  approximation of Euler-maruyama $X_{t}^{n}$ of the above mentioned equation, and we provide the rate of weak convergence $\text{Error}=\mathds{E}\left|G(X_{T})-G(X_{T}^{n})\right|$, for any function $G$ in a certain class.
\end{abstract}
\section{Introduction}
Let $X=\{X_{t}:t\geq 0\}$ be a process stochastic involving the local time defined by  the stochastic differential equations:
\begin{equation}\label{e2}
X_{t}=X_{0}+\int_{0}^{t}\varphi(X_{s})dB_{s}+\int_{\mathds{R}}\nu(da)L_{t}^{a}
\end{equation}
 where $B$ is a one-dimensional Brownian motion, $\varphi:\mathds{R}\rightarrow \mathds{R}$ is a bounded measurable function, and $\nu$ is a
 bounded measure on $\mathds{R}, L_{t}^{a}$ denotes the local time at $a$ for the time $t$ of the semimartingale X.\\
In [\ref{bib14}], if  $\nu$ is absolutely continuous with respect to the Lebesgue measure on $\mathds{R}$ (i.e.$\nu(da)=g(a)da)$
then (\ref{e2}) becomes the usual It\^o equation :
\begin{equation}\label{2}
X_{t}=X_{0}+\int_{0}^{t}\varphi(X_{s})dB_{s}+\int_{0}^{t}(g\varphi^{2})ds
\end{equation}
In general, we consider the following stochastic differential equation(SDE) with coefficients $b$ and $\sigma$, driven by
a Brownian motion $B$ in $\mathds{R}$:
\begin{equation}\label{e5}
X_{t}=X_{0}+\int_{0}^{t}\sigma(X_{s})dB_{s}+\int_{0}^{t}b(X_{s})ds
\end{equation}
where  the drift coefficient $b$ and the diffusion coefficient $\sigma$ are Borel-measurable functions from $\mathds{R}$ into $\mathds{R}$ and $X_{0}$ is an $\mathds{R}$-valued random variable, which is independent of $B$ and $b$.\\
The continuous Euler scheme $\{X^{n}_{t}, 0 \leq t\leq T \}$ for the SDE (\ref{e5}) on the time interval $[0,T]$ is defined as
follows: $X^{n}_{0}=X_{0}$, and
\begin{equation}\label{e6}
X_{t}^{n}=X_{0}^{n}+\int_{0}^{t}\sigma(X^{n}_{\eta_{k}(s)})dB_{s}+\int_{0}^{t}b(X^{n}_{\eta(k)})ds,
\end{equation}
For $\eta_{k}\leq t\leq \eta_{k+1},k=0,1,2,...,n-1$, where $0=\eta_{0}\leq \eta_{1}\leq ...\leq \eta_{n}=T$ is a sequence of
random partitions of $[0,T]$.
The weak convergence of the stochastic differential equations  has been studied by Avikainen [\ref{bib2}],
Bally and  Talay[\ref{bib3}], Mikulevi\v{c}ius and Platen[\ref{bibe16}], D.Talay and L.Tubaro[\ref{bib33}], and the weak convergence of SDEs with discontinuous coefficients has been
studied by Chan and Stramer [\ref{bibe10}], Yan [\ref{bib35}], Arturo Kohatsu-Higa [\ref{bibe21}].\\

We use the same criteria  that mentioned in [\ref{bib17}], that the definition of weak convergence with order $\gamma > 0 $ is that for all functions
$f$ in a certain class, there exists a positive constant $C$, which does not depend on $\Delta t$, such that:
\begin{equation}\label{b3}
\left|\mathds{E}[f(X_{T})]-\mathds{E}[f(X^{n}_{T})]\right|\leq C(\Delta t)^{\gamma}
\end{equation}
The stochastic differential equations of the type (\ref{e2}) have been studied by Stroock and Yor [\ref{bib21}], Portenko [\ref{bib20}], Le Gall [\ref{bib13},\ref{bib14}], Blei and Engelbert [\ref{bib7}], Bass and Chen[\ref{bib4*}] and  the approximation of Euler-Maruyama  for SDE of type (\ref{e2}) has been studied by,e.g.,  Benabdallah, Elkettani and Hiderah[\ref{bib2a}].\\
In [\ref{bib14}], when $\nu=\beta\delta_{(0)}$ ($\delta{(0)}$ denotes the Dirac measure at 0) and $\varphi=1$, we get
\begin{equation}\label{3}
X_{t}=X_{0}+B_{t}+\beta L_{t}^{0}(X),|\beta|\leq 1
\end{equation}
The solution of equation (\ref{3}) is the well-known process called the skew Brownian motion which has been studied by Harrison and Shepp [\ref{bib11}], Ouknine [\ref{bib18},\ref{bib018}], Bouhadou and Ouknine[\ref{bib}], Lejay[\ref{bib15}], Barlow [\ref{bib4}], \'Etor\'e and Martinez [\ref{bibe12}], Walsh[\ref{bib22}].\\

Our goal of this paper is that under the assumption that the SDE (\ref{e2}) has a weak solution and that it is unique, we study the conditions under which the Euler scheme
$\{X^{n}_{t}: 0 \leq t\leq T \}$ converges weakly to the exact solution $\{X_{t}:0 \leq t\leq T\}$ of the SDE (\ref{e2}).

We face two major problems: the presence of local time in equation (\ref{e2}) and the inability to provide a simple discretization scheme.
   and the presence of a discontinuity  for the diffusion coefficient $\varphi$.\\
Our paper is divided as follows: In section (\ref{sec2}), we present the  important propositions, assumptions and stochastic differential equations involving local time. Our main results on the  weak convergence of Euler scheme of one-dimensional stochastic differential equations involving the local time are given in section (\ref{sec3}). Some numerical examples are given in section (\ref{sec4}). All proofs of the theorems  are given in section (\ref{sec5}).
\section{Preliminaries and approximation}\label{sec2}
In this section, we provide the   definitions, important propositions, assumptions and stochastic differential
equations  involving local time.
\begin{defn}\label{def1}
Let $$T_{f}(x):=\sup\sum\limits_{j=1}^{N}|f(x_{j})-f(x_{j-1})|,$$
where the supremum is taken over $N$ and all partitions
$-\infty < x_{0}< x_{1}<...<x_{N}=x$  be the total variation function of $f$. Then we say that $f$ is a function of bounded
variation, if    $V(f):=\lim\limits_{x\rightarrow\infty}T_{f}(x)$ is finite, and call $V(f)$ the (total) variation of $f$.
\end{defn}
\begin{defn}
We say that a function $f$ has at most polynomial growth in $\mathds{R}$ is there exist an
integer $k$ and a constant $C >0$ such that $|f(x)| \leq C(1+|x|^k)$ for any $x \in \mathds{R}$
\end{defn}
\begin{rem}
Let $C^{k}_{p}(\mathds{R})$ denote the space of all $C^{k}$-functions of polynomial growth (together with their
derivatives).
\end{rem}
\begin{note} BV($\mathds{R}$) will denote the space of all functions $\varphi:\mathds{R}\rightarrow \mathds{R}$ of bounded variation on
$\mathds{R}$ such that:
\begin{enumerate}
  \item $\varphi$ is right continuous.
  \item There exists an $\varepsilon > 0$ such that : $\varphi(x)\geq\varepsilon$ for all $x$.
\end{enumerate}
If $\varphi$ is in $BV(\mathds{R})$, $\varphi(x^{-})$ will denote the left-limit of $\varphi$ at point $x$ and $\varphi^{'}(dx)$
 will be the bounded measure associated with $\varphi$.
\end{note}
\begin{note}
  M($\mathds{R}$) will denote the space of all bounded  measures $\nu$ on $\mathds{R}$ such that:$$|\nu(\{x\})|<1,\forall x \in \mathds{R}$$
\end{note}
We have the stochastic differential equation (\ref{e2}), if $\varphi$ is in $BV(\mathds{R})$ and $\nu$ is in $M(\mathds{R})$, then
the stochastic differential equation (\ref{e2}) has a unique strong solution as soon as $|\nu(da)| < 1 $.
\begin{thm}\label{th7}
Let $\varphi$ be in $BV(\mathds{R})$ and $\nu$ be in $M(\mathds{R})$. Then existence and pathwise uniqueness
of solution hold for (\ref{e2}).
\end{thm}
\begin{prop}\label{pr3}
Let $\nu$ be in M($\mathds{R}$). There exists a function $f$ in BV($\mathds{R}$), unique up to a multiplicative constant, such that:
\begin{equation}\label{e14}
f^{'}(dx)+(f(x)+f(x^{-}))\nu(dx)=0
\end{equation}
If we require that $f(x)\xrightarrow[x\rightarrow -\infty]{}1$, then $f$ is unique and is given by:
\begin{equation}\label{e15}
f(x)=f_{\nu}(x)=\exp(-2\nu^{c}(]-\infty ; x]))\prod\limits_{y\leq x}\left(\frac{1-\nu(\{y\})}{1+\nu(\{y\})}\right)
\end{equation}
Where $\nu^{c}$ denotes the continuous part of $\nu$.
\end{prop}
\begin{prop}\label{pr4}
Let $\varphi$ be in $BV(\mathds{R})$ and $\nu$ be in $M(\mathds{R})$ and $f_{\nu}$ be define by(\ref{e15}) and set:
\begin{equation}\label{e16}
F_{\nu}(x)=\int_{0}^{x}f_{\nu}(y)dy
\end{equation}
Then $X$ is a solution of  equation (\ref{e2}), if and only if $Y:=F_{\nu}(X)$ is a solution of:
\begin{equation}\label{e17}
dY_{t}=(\varphi f_{\nu})\circ F^{-1}(Y_{t})dB_{t}
\end{equation}
\end{prop}
The Euler scheme ${Y^{n}_{t}:0 \leq t\leq T }$ for the SDE (\ref{e17}) on the time interval $[0,T]$ is defined as
follows: $Y^{n}_{0}=Y_{0}$, and
\begin{equation}
Y^{n}_{t}=Y^{n}_{\eta_{k}}+(\varphi f_{\nu})\circ F^{-1}(Y^{n}_{\eta_{k}})(B_{t}-B_{\eta_{k}})
\end{equation}
for $\eta_{k}<t\leq \eta_{k+1}, k=0,1,2,...,n$ where $0=\eta_{0}\leq \eta_{1}\leq ...\leq \eta_{n}=T$. This Euler scheme can be written as
\begin{equation}\label{e18}
Y^{n}_{t}=Y^{n}_{0}+\int_{0}^{t}(\varphi f_{\nu})\circ F^{-1}(Y^{n}_{\eta_{k}})dB_{s}
\end{equation}
In [\ref{bib4*}], the SDE(\ref{e2}) has existence of strong solutions and pathwise uniqueness when $\varphi$ is a function on $\mathds{R}$
that is bounded above and below by positive constants and such that there is a strictly increasing function $f$ on $\mathds{R}$ such that
$$\left|\varphi(x)-\varphi(y)\right|^{2}\leq \left|f(x)-f(y)\right|,x,y \in \mathds{R}$$
  and $\nu$ is a finite measure with $|\nu(\{a\})|\leq 1$ for every $a\in  \mathds{R}$.
\begin{thm}\label{thbos}
    Suppose $\varphi$ is a measurable function on $\mathds{R}$ that is bounded
above and below by positive constants and suppose that there is a strictly
increasing function f on R such that
\begin{equation}\label{boss1}
\left|\varphi(x)-\varphi(y)\right|^{2}\leq \left|f(x)-f(y)\right|,x,y \in \mathds{R}
\end{equation}
For any finite signed measure $\nu$ on $\mathds{R}$ such that $\nu(\{x\})< \frac{1}{2}$ for each $x \in \mathds{R}$
and every$ x_{0}\in \mathds{R}$, the SDE
\begin{equation}\label{boss2}
X_{t}=X_{0}+\int_{0}^{t}\varphi(X_{s})dB_{s}+\int_{\mathds{R}}\nu(da)L_{t}^{a}, t\geq 0
\end{equation}
has a continuous strong solution and the continuous solution is pathwise unique.
  \end{thm}
  Define

Let
$$\pi(x):= \begin{cases}
      -\frac{\log(1-2x)}{2x}, & \mbox{if } x\in (-\infty,0)\cup(0,\frac{1}{2}) \\
      1, & \mbox{if } x=0.
    \end{cases}$$
    Let $\mu(dx):=\pi(\nu(\{x\})\nu(dx)$ which is a finite signed measure. Define
    $$S(x):=\int_{0}^{x}e^{-2\mu(-\infty,y]}dy$$
Since $\mu$ is a finite measure, $S^{'}$ is right continuous and strictly positive. Hence
$S$ is increasing and one-to-one. Let $S^{-1}$ denote the inverse of $S$, let $S^{'}_{\ell}$ denote
the left continuous version of $S^{'}$, i.e., the left hand derivative of s, and let $\sigma^{'}_{\ell}$ denote the left hand derivative
of $S^{-1}$. Since $\mu$ is a finite signed measure, $S^{'}$ is
of bounded variation. Let $\{Y_{t},t\leq 0\}$ solve
\begin{equation}\label{boss3}
  dY_{t}=(S^{'}_{\ell}\varphi)\circ S^{-1}(Y_{t})dB_{t},Y_{0}=S(X_{0})
\end{equation}
Let $X=S^{-1}(Y)$,  and $Y_{t}$ define by (\ref{boss3}), in [\ref{bib4*}], they must show that $X$ is a solution to (\ref{boss2}).
\section{The main results}\label{sec3}
In this section, we provide the main theorems. we study two cases: the continuous of function $\left(\varphi.f_{\nu}\right)\circ F^{-1}(.)$ and the discontinuous of function $\left(\varphi.f_{\nu}\right)\circ F^{-1}(.)$.
\subsection{Main Theorems}
In this section, we present the following results on the rates of the Euler-Maruyama approximation.\\
Let $X_{t}$ be defined as in equation (\ref{e2}), $X_{t}^{n}$ be the Euler scheme for equation (\ref{e2}),  let $F$ is defined by equation (\ref{e16}), and let $f_{\nu}$ be from the proposition(\ref{pr3}).
Here, we suppose that the assumptions of  proposition(\ref{pr4}) satisfied, and the constant $C$ may change from line to line and from theorem to theorem.
\begin{thm}\label{th1}
For any function $g:\mathds{R}\rightarrow \mathds{R}$, if for each  $G:=g\circ F^{-1} \in C_{P}^{2(\gamma+1)}(\mathds{R})$,
 and $\psi(.):=\left(\varphi.f_{\nu}\right)\circ F^{-1}(.)$ is continuous, then,
there exists a constant $C>0$, which does not depend on $n$, such that
$$ \mathds{E}\left|g(X_{T}^{n})-g(X_{T})\right|\leq \frac{C}{n^{\gamma}}$$
holds.
\end{thm}
Next we give a sufficient condition under which the Euler scheme  converges weakly to the weak solution of SDE (\ref{2}) case: discontinuous of function $\left(\varphi.f_{\nu}\right)\circ F^{-1}(.)$.
Let $\psi_{1}(z):=\liminf\limits_{y\rightarrow z}\psi^{2}(y)>0$
\begin{thm}\label{th3}
  Suppose that $\psi:=(\varphi.f)\circ F^{-1}(.)$ has at most linear growth with $D_{\psi}$ of Lebesgue measure zero and $\mathds{E}(Y_{0})^{4}<\infty$.
  If $\psi_{1}(z)>0$ for $z\in D_{\psi}$,and $G:=g\circ F^{-1} \in C_{P}^{2(\gamma+1)}(\mathds{R})$, then the Euler scheme of SDE(\ref{e2}) converges weakly to the unique weak solution of SDE(\ref{e2}), where $D_{\psi}$ is the set of discontinuous points of $\psi(y)$.
\end{thm}
If we suppose that the conditions of  theorem (\ref{thbos}) satisfied, we show the following note.
\begin{note}
The same results in the theorems (\ref{th1},\ref{th3}) stay valid, if the conditions of theorem (\ref{thbos}) are hold.
\end{note}
\section{Some examples}\label{sec4}
\begin{exmp}
In equation (\ref{e2}), let  $\nu(dx)=\alpha \delta_{0}(dx), |\alpha|<1$  where ($\delta_{0}$ is the Dirac measure  at 0 ) and
  \begin{equation*}
  \varphi=\begin{cases}
   \frac{1+\alpha}{1-\alpha}\exp \left(-\frac{1-\alpha}{1+\alpha}x\right) & \mbox{ if } x  \geq 0\\
 \exp \left(x\right)& \mbox{ if } x< 0
\end{cases}
  \end{equation*}
Using proposition (\ref{pr3}), we have
$$f_{\nu}(x)=\begin{cases}
\frac{1-\alpha}{1+\alpha} & \mbox{if  } x  \geq 0\\
 1& \mbox{if } x< 0
\end{cases}$$
Using proposition (\ref{pr4}), we have
  \begin{equation*}
F(x)=\begin{cases}
\frac{1-\alpha}{1+\alpha}x & \mbox{if  }  x  \geq 0\\
 x& \mbox{if }  x< 0
\end{cases}
\end{equation*}
Then $$\psi(y):=(\varphi.f)\circ F^{-1}(y)=\begin{cases}
\exp \left(-y\right) & \mbox{if } y  \geq 0\\
 \exp \left(y\right)& \mbox{if } y< 0
 \end{cases}$$
 We note $\psi$ is continuous, then for any function $g:\mathds{R}\rightarrow \mathds{R}$, such that $g\circ F^{-1}\in C^{2(\gamma+1)}_{p}(\mathds{R})$, for any $\gamma >0$, for example
 $$g(x)=\frac{1}{1+\left(\frac{1-\alpha}{1+\alpha}x\right)^{2}}.1_{x\geq 0}+\frac{1}{1+x^{2}}1_{x<0}$$
  there exist a constant $C>0$, which does not depend on $n$, such that
\begin{center}
$ \left|\mathds{E}\left[g(X_{T}^{n})\right]-\mathds{E}\left[g(X_{T})\right]\right|\leq \frac{C}{n^{\gamma}}$
\end{center}
\end{exmp}
\begin{exmp}
 In equation (\ref{e2}), let $\nu(x)=\alpha\delta_{0}(x),|\alpha|<1$, where $\delta_{0}$ is the Dirac measure  at 0, and $\mathds{E}(X_{0})^{4}=\mathds{E}(Y_{0})^{4}<\infty$
we define $\varphi=1$, such that\\
Using proposition (\ref{pr3}), we have
\begin{flalign*}
f(x)=f_{\nu}(x)=\begin{cases}
\frac{1-\alpha}{1+\alpha} & \mbox{if  }  x\geq 0\\
 1 & \mbox{if  } x< 0
\end{cases}
\end{flalign*}
And by using proposition (\ref{pr4}), we have\\
\noindent\begin{tabularx}{\textwidth}{@{}XXX@{}}
  \begin{equation*}
  F_{\nu}(x)=\begin{cases}
\frac{1-\alpha}{1+\alpha}x & \mbox{if  }  x\geq 0\\
 x & \mbox{if  } x< 0
\end{cases}
\end{equation*}&
\begin{equation*}
F^{-1}_{\nu}(y)=\begin{cases}
\frac{1+\alpha}{1-\alpha}y & \mbox{if  }  y\geq 0\\
 y & \mbox{if  } y< 0
\end{cases}
\end{equation*}
\end{tabularx}
Then, we obtain
\begin{equation}
  \psi(y)=\left(\varphi.f_{\nu}\right)\circ F^{-1}(y)=\begin{cases}
\frac{1-\alpha}{1+\alpha} & \mbox{if  }  y\geq 0\\
 1 & \mbox{if  } y< 0
\end{cases}
\end{equation}
And
\begin{align}\label{e26}
Y_{t}&=Y_{0}+\int_{0}^{t}\psi(Y_{s})dB_{s}
\end{align}
The Euler scheme $\left\{Y^{n}_{t}:0 \leq t\leq T \right\}$ for the SDE (\ref{e26}), on the time interval $[0,T]$ is defined as
follows: $Y^{n}_{0}=Y_{0}$, and
\begin{align}\label{e27}
Y^{n}_{t}&=Y^{n}_{\eta_{k}}+\psi(Y^{n}_{\eta_{k}})(B_{t}-B_{\eta_{k}})
\end{align}
For $\eta_{k}<t\leq \eta_{k+1},k=0,1,2,...,n$, where $0=\eta_{0}\leq \eta_{1}\leq ...\leq \eta_{n}=T$,
the coefficient of equation (\ref{e26}) is discontinuous at point 0,
and, we have:
\begin{itemize}
  \item $D_{\psi}=\{0\}$
  \item $\psi_{1}(z):=\liminf\limits_{z\rightarrow y}\psi^{2}(y)>0$
  \item $\psi(\cdot)$ is at most linear growth with $D_{\psi}$
\end{itemize}
We have the Euler scheme of SDE(\ref{e27}) converges weakly to the unique weak solution of SDE((\ref{e26})), then
the Euler scheme $X_{t}^{n}$ converges weakly to the unique weak solution$X_{t}$ of SDE((\ref{e2})).
\end{exmp}
\section{Proofs of theorems}\label{sec5}
Before proving the theorems below, we introduce some notations.
\begin{note}
  The constant $C$ may change from line to line and from theorem to theorem.
\end{note}
Let $B$ be a one-dimensional brownian motion
$(B_{t})_{t\in[0,T]}$, and $\varphi$ is in $BV(\mathds{R})$ and $\nu$ is in $M(\mathds{R})$ and
\begin{equation}\label{bbe1}
X_{t}=X_{0}+\int_{0}^{t}\varphi(X_{s})dB_{s}+\int_{\mathds{R}}\nu(da)L_{t}^{a}.
\end{equation}
and $Y_{t}$ is solution of the equation
\begin{equation}\label{bbe2}
dY_{t}=\left(\varphi.f_{\nu}\right)\circ F^{-1}(Y_{t})dB_{t}.
\end{equation}
for which, as we get by the proposition(\ref{pr4}), where $F$ be defined by  the equation (\ref{e16}), $f_{\nu}$ get by the proposition(\ref{pr3}),
and let $Y_{t}^{n}$ is a solution by Euler scheme for equation(\ref{e17}) and defined by equation(\ref{e18}).
\subsection{Proof of theorem (\ref{th1})}
\begin{defn}\label{def2}
In [\ref{bib17}], a time discrete approximation $X^{n}$  converges weakly with order $\gamma>0$ to $X$ at time $T$ as $n\rightarrow \infty$,
if for each  $G\in C_{P}^{2(\gamma+1)}(\mathds{R})$,
there exists a constant $C>0$, which does not depend on $n$, such that
$$\left| \mathds{E}(G(X_{T}^{n}))-\mathds{E}(G(X_{T}))\right|\leq \frac{C}{n^{\gamma}} \text{holds.}$$
\end{defn}
Using the above definition, we can prove theorem (\ref{th1}).\\
\begin{enumerate}
                   \item The proposition (\ref{pr4}) is satisfied , then $X_{t}=F^{-1}(Y_{t})$ is solution uniqueness of the SDE (\ref{e2}).
                   \item $G:=g\circ F^{-1}\in C^{2(\gamma+1)}_{p}$
                   \item By using definition (\ref{def2}), we have\\
                   \begin{align*}
                     \left|\mathds{E}[g(X_{T})]-\mathds{E}[g(X^{n}_{T})]\right|&=
                   \left|\mathds{E}[g\circ F^{-1}(Y_{T})]-\mathds{E}[g\circ F^{-1}(Y^{n}_{T})]\right|\\
                      & =\left|\mathds{E}[G(Y_{T})]-\mathds{E}[G(Y^{n}_{T})]\right|\leq \frac{C}{n^{\gamma}}
                   \end{align*}
                 \end{enumerate}
\subsection{Proof of theorem (\ref{th3})}
If we have SDE of type \begin{equation}\label{b}
Y_{t}=Y_{0}+\int_{0}^{t}\sigma(Y_{s})dB_{s}
 \end{equation}
 and the Euler scheme of equation (\ref{b}) give by
\begin{equation}\label{b1}
  Y_{t}=Y_{0}+\int_{0}^{t}\sigma(Y_{\eta_{n}(s)})dB_{s}
\end{equation}
Here, we present the following theorem(see [\ref{bib35}]):
\begin{thm}\label{th4}
  Suppose that $\sigma(y)$ has at most linear growth with $D_{\sigma}$ of Lebesgue measure zero and $\mathds{E}(Y_{0})^{4}<\infty$.
  If $\sigma_{1}(z):=\liminf\limits_{z\rightarrow y}\sigma^{2}(y)>0$ for $z\in D_{\psi}$, then the Euler scheme of SDE(\ref{b}) converges weakly to the unique weak solution of SDE(\ref{b}), where $D_{\sigma}$ is the set of discontinuous points of $\sigma(y)$.
\end{thm}
Using the above theorem, we can prove theorem (\ref{th3}).\\
\begin{enumerate}
                   \item The proposition (\ref{pr4}) is satisfied , then $X_{t}=F^{-1}(Y_{t})$ is solution uniqueness of the SDE (\ref{e2}).
                   \item $G:=g\circ F^{-1}\in C^{2(\gamma+1)}_{p}$
                   \item $(\varphi.f_{\nu})\circ F^{-1}(\cdot)$ has at most linear growth with $D_{\sigma}$ of Lebesgue measure zero.
                   \item $\mathds{E}(Y_{0})^{4}<\infty$
                   \item By using definition (\ref{def2}) and theorem (\ref{th4}) , we have\\
                    \begin{align*}
                   \left|\mathds{E}[g(X_{T})]-\mathds{E}[g(X^{n}_{T})]\right|&=
                   \left|\mathds{E}[g\circ F^{-1}(Y_{T})]-\mathds{E}[g\circ F^{-1}(Y^{n}_{T})]\right|\\
                   &=\left|\mathds{E}[G(Y_{T})]-\mathds{E}[G(Y^{n}_{T})]\right|\leq \frac{C}{n^{\gamma}}
                   \end{align*}
                 \end{enumerate}
\bibliographystyle{plain}

\end{document}